\documentstyle{article}

 \def\vt{t\kern-0.22em\raise.18ex\hbox{\char'47}\lower.18ex\hbox{}\kern-0.08em}
 \newtheorem{pr}{Proposition}[section]
\newtheorem{th}{Theorem}[section]

\newtheorem{ex}{Example}[section]

\newcommand{\old}[1]{{}} 

\newcounter{obr}
\newcounter{tabul}
\begin{document}
\title{
Clique Cover Width and Clique Sum\footnote{This paper appeared in Congressus Numerantium, 218 (2013), 135-140.}}
\author{Farhad Shahrokhi\\
Department of Computer Science and Engineering,   UNT\\
farhad@cs.unt.edu
}

\date{}
\maketitle
\thispagestyle{empty}
\date{} \maketitle

\begin{abstract}
For a clique cover $C$ in the undirected graph $G$, the clique cover
graph of $C$ is the graph obtained by contracting the vertices of each
clique in $C$ into a single vertex. The clique cover width of G, denoted
by $CCW(G)$, is the minimum value of the bandwidth of all clique
cover graphs of $G$. When $G$ is the clique sum of $G_1$ and $G_2$, we prove
that $CCW(G) \le 3/2(CCW(G_1) + CCW(G_2))$.
\end{abstract}

\section{Introduction and Summary}
Throughout this paper, $G = (V (G),E(G))$ denotes a graph on $n$ vertices.
$G$ is assumed to be undirected, unless stated otherwise. Let 
$L =
\{v_0, v_1, ..., v_{n-1}\}$ be a linear ordering of vertices in $G$. The width of $L$, denoted
by $w(L)$, is $max_{v_i,v_j\in E(G)}\{|j-i|\}$. The bandwidth of G, denoted by
$BW(G)$, is the smallest width of all linear orderings of $V(G)$ [1], [2], [3].
Unfortunately, computing the bandwidth is NP −hard, even if $G$ is a tree
[6]. A clique cover $C$ in G is a partition of $V(G)$ into cliques. Throughout
this paper, we will write $C = \{c_0, c_1, ..., c_t\}$ to indicate that $C$ is an ordered
set of cliques. For a clique cover $C = \{c_0, c_1, ..., c_t\}$, in $G$, let the width of $C$,
denoted by $w(C)$, denote 
$max\{|j-i||xy\in  E(G), x\in  c_i, y\in  c_j, c_i, c_j\in  C\}$.
The clique cover width of G denoted by $CCW(G)$, was defined by the author
in [10], to be the smallest width of all ordered clique covers in $G$. Equivalently,
one can define $CCW(G)$ to be the minimum value of the bandwidth
overall clique cover graphs of G, as stated in the abstract. Combinatorial
properties of CCW(G) have been studied. Clearly, $CCW(G)\le BW(G)$.
In addition, it is a easy to verify that $BW(G)\le \omega(G).CCW(G)$, where
$\omega(G)$ is the size of a largest clique in $G$, and that any G with $CCW(G) = 1$.
is an incomparability graph [10]. Furthermore, we proved in [11] that
$CCW(G)\ge \lceil s(G)/2\rceil -1$ 
 for any graph $G$, where $s(G)$ is the largest number
of leaves in an induced star in $G$, and that for any incomparability graph
$G$, 
$CCW(G)\le s(G)-1$. Additionally, in [11], we further explored the
class of unit incomparability graphs, or graphs $G$ with $CCW(G)=1$; this
class is fairly large and contains the classes of the unit interval graphs and
co-bipartite graphs. Furthermore, in [11] we introduced the unit incomparability
dimension of G, or Udim(G), which is a parameter similar to the
cubicity of G [7], [14], [15]. Specifically, Udim(G) is defined as the smallest
integer d so that G is the intersection graph of d unit incomparability
graphs. In [11], we also proved a decomposition theorem that establishes
the inequality $Udim(G)\le CCW(G)$, for any graph G. The upper bound
is improved from $CCW(G)$ to $O(log(CCW(G)))$ in [13]. Finally, we have just
proved that every planar G is the intersection graph of a chordal graph and
a graph whose clique cover width is at most seven [12].
The main application of the clique cover width is in the derivation of
separation theorems in graphs; particularly the separation can be defined
for more general types measures [9], instead of just the number of vertices.
For instance, given a clique cover $C$ in $G$, can $G$ be separated by removing a
*small* number of cliques in $C$ so that each the two remaining subgraph of
G can be covered by at most $k|C|$ cliques from $C$, where $k< 1$ is a constant
[9, 13]? Our recent work shows a close connection between the tree width
of $G$, or $tw(G)$, and $CCW(G)$. Recall that $tw(G)-1$ is the minimum of
the maximum clique sizes of all chordal graphs that are obtained by adding
edges to $G$ [8],[4].
Let $G_1$ and $G_2$ be graphs so that $V(G_1)\cap V(G_2)$ is a clique in both
$G_1$ and $G_2$. Then, the clique sum of $G_1$ and $G_2$, denoted by $G_1\oplus G_2$, 
is a graph $G$ with $V(G) = V(G_1)\cap V(G_2)$, and $E(G) = E(G_1)\cup E(G_2)$.
Clique sums are intimately related to the concept of the tree width and tree
decomposition; Specifically, it is known that if the tree widths of $G_1$ and 
$G_2$
are at most $k$, then, so is the tree width of $G_1\oplus G_2$ [5]. Unfortunately, this
is not true for the clique cover width. As seen in the following example.

\begin{ex}\label{e1}
{\sl 
 Let $P_1$ and $P_2$ be paths on $2t + 1, t\ge 1$ vertices, then,
$CCW(P1)=CCW(P2) = 1$. Now select the unique vertex $x$ in the middle
of the two paths and take the sum of $P_1$ and $P_2$ at $x$. Then, 
$s(P_1\oplus P_2) = 4$,
and consequently, $CCW(P_1\oplus P_2)\ge s(G)/2 = 2$, whereas, a simple ordering
linear ordering of vertices shows that lower bound is achievable, and in fact
$CCW(P1\oplus P2)=2$.
}
\end{ex}

In this paper we study the clique cover width of the clique sum of two
graphs and establish the inequality $CCW(G_1\oplus G_2)\le (3/2)(CCW(G_1) +
CCW(G_2))$.
In section two we define some Specific technical concepts that will be
used to

\section{Preliminaries}

Let $C = \{c_0, c_1, ..., c_t\}$ be a clique cover in $G$. Any set of consecutive cliques
in $C$ is called a strip. 
Any strip of of cardinality  $w=w(C)$ is called a {\it block}. Let $S$ be a strip  in $C$. We denote
by $C^l(S)$ and $C^r(S)$, the largest strips that are entirely to the left, and to
the right of $S$, respectively; thus, $C=\{C^l(S),S,C^r(S)\}$. Note that if $S$ is a block, then,
the removal of cliques in $S$ disconnects the subgraphs induced on 
$C^l(S)$ and $C^r(S)$ in $G$. For a strip  $S$ in $C$, we denote by 
$S^l$  and $S^r$ the largest
strips of cardinality at most $w$ that immediately precede and proceed
$S$, in $C$, respectively.
 Let $B=\{c_k,c_{k+1},...,c_{k+w-1}\}, k\ge 0$   be a block in a clique cover 
$C =
\{c_0, c_1, ..., c_t\}$ for $G$. We will define a partition of $C$, denoted by 
$P(C,B)$, as follows. Let $C^l(B) =\{c_0, c_1, ..., c_{k-1}\}$, and 
$C^r(B) =\{c_{k+w}, c_{k+w+1}, ..., c_t\}$. Now, let $k=p.w+r, r\le w-1$
  Define $S_0$ to be the set of first
$r$  consecutive cliques in 
$C^l(B)$. For $i = 1, 2, ..., p$, let 
$S_i$  be the block
in $C$ starting at $c_i{(i−1).w+r}$. It is easy to verify that 
$\{S_0, S_1, ..., S_p\}$ is a
partition of $C^l(B)$. Similarly, one can construct a partition of $C^r(B)$ 
of form
$\{S_{p+2}, S_{p+3}, ..., S_q\}$, for some properly defined 
$q$, so that $S_i$ is a block for 
$i = p+2, p+3, ..., q$, and the width of $S_p$  is at most 
$w(C)$ . By combining these
two partitions with $S_{p+1} = B$, as the middle part, we obtain the partition
$P(C,B)$ of $C$ into strips. For $i, j = 0, 1, ..., q$, we define the distance of
strips $S_i$ and $S_j$ , in the partition $P(C,B)$, to be $|j-i|$.

\begin{pr}\label{p1}
{\sl
 Let C be a clique cover in $G$ and let $B$ be a block in $C$.
Then, there is an ordered partition $P(C,B) = {S_0, S_1, ..., S_p,B, S_p+2, ..., S_q}$
of $C$ into strips so that $\{S_0, S_1, ..., S_p\}$ partitions $C^l(B)$, 
$B = S_p+1$, and
${S_p+2, ..., S_q}$ partitions $C^r(B)$. Moreover, the first and last strips in $P(C,B)$
have widths at most $w(C)$, whereas, the remaining elements of $P(C,B)$ are
blocks.
}
\end{pr}

\section{Main Results}
Let $G=G_1\oplus G_2$, and let $C_1$ and  $C_2$ be clique covers in $G_1$ and $G_2$, respectively. 
Let $S_1$ and $S_2$ be strips in $C_1$ and $C_2$, respectively.
The interleave ordering of $S_1$ and $S_2$, denoted by,
$S_1\oplus  S_2$ is an ordered set of cliques obtained by placing the 
first clique in $S_1$ after the first clique in $S_2$, the second clique in 
$S_1$ after the second clique in $S_2$, etc., until all cliques in one $S_i, i = 1, 2$, are used, then, one places all
remaining cliques in $S_{i+1} (mod 2)$ in $S_1\oplus  S_2$. 
Note that $S_1\oplus \emptyset= S_1$, 
and $\emptyset\oplus  S_2 = S_2$. 

\begin{th}\label{t1}
{\sl 
 Let $G_1$ and $G_2$ be graphs, where $S= V(G_1)\cap  V(G_2)$ 
induces a clique in $G_1$ and $G_2$. Then, $CCW(G_1\oplus  G_2)\le 3/2(CCW(G_1) +
CCW(G_2))$.
}
\end{th}

{\bf Proof}. To prove the claim, let $C_1 = \{c_0, c_1, ..., c_a\}$ and 
$C_2 = \{c'_0, c'_1, ..., c'_b\}$
be clique covers in $G_1$  and $G_2$, respectively; it suffices to construct a clique
cover $C$ for $G$, of width at most $3/2(w(C_1) +w(C_2))$.
. Let $B_1$ and $B_2$ be the blocks in
$C_1$ and $C_2$ that contain all vertices in $S$, respectively. Note that such blocks
must exist, since $S$ induces a clique in $G_1$, as well as, in $G_2$. 
Construct the
partitions $P(B_1,C_1) = \{S_0, S_1, ..., S_p,B_1, S_{p+2}, ..., S_q\}$ and 
$P(B_2,C_2)=
\{T_0, T_1, ..., T_r,B_2, T_{r+2}, ..., T_u\}$. 
Now, interleave $B_1$ and $B_2$, to obtain a
strip $I(B_1,B_2)$. Next, for 
$i = p, p − 1, ..., 0, j=r, r − 1, ..., 0$, if 
$p − i=
r − j$, then, interleave $S_i$ and $T_j$ , to obtain the strip 
$I(Si, Tj)$. (So, we
interleave those strips in $c_1$ and $C2$ that are to left of $B_1$ and $B_2$, and
are of the same distance from $B_1$ and $B_2$, respectively.) 
Let $L^l$ denote the  
ordered union of all these nterleaved strips.  
Similarly, for
$i = p+2, p+3, ..., S_q, j = r +2, r +3, ..., u$, if 
$p−i = r −j$, then, interleave
$S_i$ and $T_j$ to obtain $I(S_i, T_j)$ , and, let  
$L^r$ denote the  
ordered union of all these nterleaved strips.  
Now, let $L=\{L^l,I(B_1,B_2),L^r\}$.
Note that $L = \{l_0,l_1, ..., l_{a+b+1}\}$ is an ordered clique cover for $G$, 
but the
cliques in $L$ are not disjoint. It is important to observe that any clique in
$L$ either belongs to $C_1$ or to $C_2$, and that the relative orderings of cliques
in $C_1$ or $C2$ remain the same in $L$.

{\bf Claim.}  Let $xy\in E(G_1)\cup  E(G_2)$, 
with $x\in l_r, y\in l_t, l_r,l_t\in L$.
Then, $|t- r|\le   w(C_1) + w(C_2)-1$.

{\bf Proof.} Clearly, $l_r, l_t\in   C_i$, for some $i = 1, 2$,  If $l_r, l_t$  are in the same strip
of $L$, as prescribed in the interleave construction, then, the claim is true.
So assume that $l_r$ and $l_t$ are in different strips of $L$. To verify the claim
it suffices to verify that $l_r$ and $l_t$ can only be in adjacent strips. However,
this is a consequence of the interleave construction. $\Box$ 

We will now convert $L$ to a clique cover in G as follows. Create a new
clique in the middle of the strip $I(B_1,B_2)$, place $S$ in this clique, and remove
all vertices of $S$ from any cliques in $L$. Let $C$ denote the clique cover that
	is constructed this way. Using the above claim, we get, $w(C)\le  w(C_1) +
w(C_2)-1 + {w(C_1)+w(C_2)\over 2}  + 1$. Consequently, $w(C)\le   3/2(w(C_1) + w(C_2))$,
as stated.$\Box$


\begin{thebibliography}{99}

\bibitem{1}  J. B¨ottcher, K. P. Pruessmannb, A. Taraz, A. W¨urfel, Bandwidth,
treewidth, separators, expansion, and universality, European Journal
of Combinatorics, 31, 2010, 1217, 1227.

[2] P.Z. Chinn, J. Chvatalov , A. K. Dewdney, N.E. Gibbs, Bandwidth
problem for graphs and matrices- A Survey, Journal of Graph Theory,
6(3), 2006, 223-254

[3] J. Diaz, J. Petit, M. Serna A survey of graph layout problems, ACM
Computing Surveys (CSUR) 34(3), 2002, 313 - 356.

[4] H.L. Bodlaender , A Tourist Guide through Treewidth. Acta Cybern.
11(1-2), 1993, 1-22.

[5] H.L. Bodlaender, A partial k-arboretum of graphs with bounded
treewidth. Theoretical Computer Science, 1998.

[6] M.R. Garey and D.J. Johnson, Computers and Intractability: A Guide
to the Theory of NP-Completeness, Freeman, San Francisco, CA, 1978

[7] F. Roberts. Recent Progresses in Combinatorics, chapter On the boxicity
and cubicity of a graph, pages 301-310. Academic Press, New
York, 1969.

[8] N. Robertson, P. D. Seymour, Graph minors III: Planar tree-width,
Journal of Combinatorial Theory, Series B *36*, 1984, (1): 49-64.

[9] F. Shahrokhi, A New Separation Theorem with Geometric Applications,
Proceedings of EuroCG2010, 2010, 253-256.

[10] F. Shahrokhi, On the clique cover width problem, Congressus Numerantium,
205 (2010), 97-103.

[11] F. Shahrokhi, Unit Incomparability Dimension and Clique Cover
Width in Graphs, Congressus Numerantium, 213 (2012), 91-98.

[12] F. Shahrokhi, New representation results for planar graphs, 29th European
Workshop on Computational Geom., 2013.

[13] F. Shahrokhi, in preparation.

[14] L. Sunil Chandran and Naveen Sivadasan. Boxicity and treewidth.
Journal of Combinatorial Theory, Series B, 97(5):733-744, September
2007.

[15] Abhijin Adiga, L. Sunil Chandran, Rogers Mathew: Cubicity, Degeneracy,
and Crossing Number. FSTTCS 2011: 176-190.

\end{thebibliography}
\end{document}